# Existence of Least-perimeter Partitions

Frank Morgan[*]
Department of Mathematics and Statistics
Williams College
Williamstown, MA 01267, USA

## Abstract

We prove the existence of a perimeter-minimizing partition of $\mathbf{R}^n$ into regions of unit volume.

## 1. Introduction

In 1894 Lord Kelvin ([1], [2, Chapt. 15]) conjectured that certain relaxed truncated octahedra provide the least-perimeter way to partition $\mathbf{R}^3$ into unit volumes. One hundred years later Weaire and Phelan ([3], [2, Chapt. 15]) found a better candidate, a mixture of relaxed dodecahedra and 14-hedra, now conjectured to be the minimizer. There is no proof in sight. It has been an open question whether a minimizer exists. Since space is infinite, one has to consider a limiting perimeter-to-volume ratio, and all estimates could disappear to infinity in a purported optimal limit configuration. Worse, in theory each region of unit volume could have many far-flung components, some of which could disappear to infinity, resulting in regions of volume less than one. First we show how to construct a configuration which has optimal limiting ratio although some regions might have volume less than one. Second we show how to readjust to unit volumes by identifying pieces of some regions with others. Our construction yields regions which are bounded but not connected.

In $\mathbf{R}^2$, regular hexagons provide a perimeter-minimizing partition of the plane into unit areas, as was proved by Hales ([4], [2, Chapt. 15]) in 1999. The solution is not unique. Indeed, compact alterations do not affect the limiting perimeter-to-volume ratio and yield other minimizers, some of which have disconnected regions.

For perimeter-minimizing clusters of *finitely many* regions of prescribed volumes (see [2, Chapt. 13]), existence and partial regularity were proved by Almgren [5]. In theory one has to allow the possibility that a region may be disconnected, because distinct components may be connected by thin tubes of arbitrarily small perimeter. This possibility of regions of several components was the main difficulty in the proof of the Double Bubble Conjecture ([6–9]), which says that the familiar double soap bubble is the least-perimeter way to enclose and separate two prescribed volumes, newly extended by Reichardt [9] to $\mathbf{R}^n$.

[*]Frank.Morgan@williams.edu

## 2. Existence of Perimeter-minimizing Partitions

**Theorem.** *There exists a partition of $\mathbf{R}^n$ into regions of unit volume minimizing the limit superior of the perimeter per volume ratio for balls about the origin.*

By "perimeter" we mean the $(n-1)$-dimensional Hausdorff measure [2, §2.3] of the union of the topological boundaries of the regions. One considers the restrictions of the regions to a closed ball of radius r, computes the Hausdorff measure of the union of their topological boundaries, divides by the volume of the ball, and takes the limit superior as r approaches infinity. Note that the area of the truncating sphere of radius r is negligible in the limit. Alternatively, one could work in the category of the integral currents of geometric measure theory [2].

*Proof.* Choose units for perimeter such that the greatest lower bound of all such limits superior is 1. Let $\varepsilon_i = 2^{-i}$. For $i \geq 0$, choose partitions $P_i$ of $\mathbf{R}^n$ into unit volumes such that the perimeter per volume ratio $a_i(r)$ inside a closed ball of radius r satisfies

$$\limsup a_i(r) < 1 + \varepsilon_i.$$

For a rapidly increasing sequence $r_i$ with $r_0 = 0$ to be specified, consider the partition P equal to $P_i$ for $r_i < r < r_{i+1}$ together with the spheres of radius $r_i$, with perimeter per volume ratio $a(r)$ ignoring the spheres, which do not affect the lim sup. Note that P is a partition of $\mathbf{R}^n$ into regions of volumes at most 1; we'll take care of that deficiency later.

To make the final estimates come out conveniently, we'll choose the sequence $r_i$ to increase by factors of at least 5. Choose $r_1 \geq 5$ such that

(1)  $\qquad\qquad\qquad a_0(r_1) < 1 + 4\varepsilon_1.$

Choose $r_i \geq 5r_{i-1}$ such that

(2)  $\qquad\qquad\qquad a_i(r) < 1 + \varepsilon_i \qquad$ for $r \geq r_i$ ,

which holds for all sufficiently large $r_i$ , and

(3)  $\qquad\qquad\qquad a_i(r_i) > 1 - \varepsilon_i$ ,

which holds for arbitrarily large $r_i$. We will prove by induction that for $i \geq 1$,

(4)  $\qquad\qquad\qquad a(r_i) < 1 + 4\varepsilon_i$ ,

For $i = 1$, (4) holds by (1). Assuming (4) for $a(r_i)$, to estimate $a(r_{i+1})$, note that since $r_{i+1} \geq 5r_i$, on at least 4/5 of the ball of radius $r_{i+1}$, P coincides with $P_i$, which has $a_i(r_{i+1}) < 1 + \varepsilon_i$ by (2). On at most 1/5 of the ball, $P_i$ is replaced by P, which has by induction $a(r_i) < 1 + 4\varepsilon_i$, while $a_i(r_i) > 1 - \varepsilon_i$ by (3); this replacement increases the perimeter per volume ratio by less than $(1/5)(5\varepsilon_i) = \varepsilon_i$. Therefore $a(r_{i+1}) < (1 + \varepsilon_i) + \varepsilon_i = 1 + 4\varepsilon_{i+1}$.

Similarly since for $r_i \le r \le r_{i+1}$ P consists of $P_i$ replaced by P in the ball of radius $r_i$, therefore by (3), (4), and (2)

$$(5) \qquad\qquad a(r) < \frac{1 + 4\varepsilon_i}{1 - \varepsilon_i}(1+\varepsilon_i) \ .$$

By (5), lim sup $a(r) \le 1$.

Finally we show how to change P to a partition into regions of volumes exactly 1, without altering lim sup $a(r)$. Inside annular layers $\{r_i \le r \le r_{i+1}\}$, by identifying regions we may assume that at most one region in each layer has volume less than 1/2. Hence each layer has finitely many regions. Order the regions, starting with the innermost annular layer and moving outward: $R_1$, $R_2$, $R_3$, ... . Choose two regions $T_i$, $T_i'$ of volume at least 1/2 in each layer starting with the second layer; these relatively scarce $T_i$, $T_i'$ will have parts cut off and made part of the $R_i$ to restore the volume of each $R_i$ to 1. (Since each $T_i$, $T_i'$ is an $R_j$, it will eventually have its volume restored to 1 as well.)

If vol($R_1$) < 1, divide off a piece of $T_1$ of volume 1−vol($R_1$) using a portion of a sphere about the origin and identify it with $R_1$ to make its volume 1; or if necessary take all of $T_1$ and a piece of $T_1'$. Similarly take from each $T_i$, $T_i'$ pair volume 1−vol($R_i$) for $R_i$ to make its volume 1. Now all the regions have volume 1 as desired. Since the $T_i$ are in successive annuli of rapidly growing radii and the added perimeter comes from a portion of a sphere, the limiting area per volume ratio of P, lim sup $a(r)$, remains 1. Thus P is the desired partition of $\mathbf{R}^n$ into regions of unit volume minimizing the limit superior of the perimeter per volume ratio.

*Remark.* Since each region consists of a piece from one annular layer plus a piece from a pair in another annular layer, each region is bounded.

## 3. Manuel A. Fortes

It was my good fortune as a Johnny-come-lately to have the last chance to be one of Fortes's many collaborators. His kind invitation brought me on my first visit to Lisbon in September 2006. Though in his last months, Fortes was still in rare form, leading one of his famous expeditions to lunch, solicitously making sure that I had a chance to sample the best of the local cuisine, insistently raising his favorite mathematical problems and questions for me, the mathematician. For example, he conjectured a generalization to 3D of the celebrated four-color theorem, which says that a planar map can be colored with four colors, with no adjacent countries the same color. Somehow neither I nor any of my mathematician friends had ever considered this natural and compelling question. On further research, I found that there was indeed a conjecture that for convex polyhedra in 3D, six colors suffice. Sadly, it is apparently false. A recent posting to the web [10] says that no finite number of colors suffices in general.

Fortes was excited about an apparently new unstable periodic chain bubble cluster, which appeared in one of his last publications with Vaz [11].

Our most fruitful discussions concerned the pressure in planar clusters of many unit-area bubbles. The whole cluster as well as the individual bubbles tends to be hexagonal. We soon realized that the pressure is largest not in bubbles at the center of the cluster, as one might expect, but in bubbles at the corners of the cluster. The smallest pressure bubble is also on the outer boundary, halfway between the corners. The pressure at the center approaches an interesting limiting value. Fortes urged me to write a joint paper promptly. I remember his joy when I produced a manuscript that very afternoon before my departure. We submitted the proof corrections to *Philosophical Magazine Letters* the following April 18; Fortes passed away four days later. We dedicated the paper [12] to him.

## Acknowledgements

This paper was inspired by discussions in my Fall 2007 senior seminar at Williams College on "The Big Questions," including packing and partitioning problems as well as many of the Clay Institute Million-Dollar Millennium Problems. I want to thank my students: Jonathan Berch, Trubee Davison, Christophe Dorsey, Kimberly Elicker, Brian Hwang, Corey Levin, Anthony Marcuccio, Edward Newkirk, Rahul Shah, Sebastian Shterental, Sara Siegmann, Amy Steele, and Kristin Sundet.

49Q10 Optimization of shapes other than minimal surfaces